# High-accuracy finite-difference schemes for solving elastodynamic problems in curvilinear coordinates within multiblock approach

Leonid Dovgilovich[1,2], Ivan Sofronov[1,2]

[1]*Schlumberger, Moscow, Pudovkina 13,*

[2]*MIPT, Moscow region, Dolgoprudny, Institutskii per. 1*

ldovgilovich@slb.com, isofronov@slb.com

**Abstract**

We propose highly accurate finite-difference schemes for simulating wave propagation problems described by linear second-order hyperbolic equations. The schemes are based on the summation by parts (SBP) approach modified for applications with violation of input data smoothness. In particular, we derive and implement stable schemes for solving elastodynamic anisotropic problems described by the Navier wave equation in complex geometry. To enhance potential of the method, we use a general type of coordinate transformation and multiblock grids. We also show that the conventional spectral element method (SEM) can be treated as the multiblock finite-difference method whose blocks are the SEM cells with SBP operators on GLL grid.

# 1 Introduction

The summation-by-parts simultaneous-approximation-terms (SBP-SAT) approach [12], [3], [17], [18] allows generation of high-order accurate finite-difference schemes for simulating wave propagation problems so that the matrices of spatial operator for second-order hyperbolic equations are symmetric. The symmetry ensures stability of derived difference schemes under Courant-Friedrichs-Lewy (CFL) condition. Extension of the application field of such schemes is possible by using coordinate transformations that enable consideration of curvilinear grids adapted to complex geometry [1] as well as by using a multiblock approach for handling even more complex geometry and discontinuous coefficients of governing equations [11].

In this paper, we use the SBP idea to obtain some new difference schemes of integrating the scalar wave equation and anisotropic Navier wave equations in curvilinear coordinates including a multiblock framework. The original differential equations are considered in the conservative (divergence) form. Note the following specific features of proposed schemes.

1) Our difference schemes are derived on *conventional grids*, i.e., when all components of the solution vector and coefficients of governing equations belong to a single grid point (this is in contrast to *staggered grids* in which all these variables are distributed on a grid cluster whose points are shifted by the half-spacing in space and time from the physical point of medium, see [28], [23], [15]). The choice to apply conventional grids was made because, in our opinion, they are algorithmically easier for using curvilinear coordinates and a multiblock approach. In addition, the SBP approach is traditionally developed on conventional grids.

2) Unlike methods [1], [16], and [18], for the second-order equations we use only operators of the first derivatives while discretizing all other derivatives. One advantage of this approach consists of the possibility of implementing the scheme with an optimal number of the operations per grid point in the calculation of the mixed derivatives: $O(p)$ where $p$ is the number of points in one direction of the stencil (instead of $O(p^2)$ as it would be for approximation of mixed derivatives directly by second-order difference operators).

3) As is known, the disadvantage of using central-difference operators of the first derivatives for constructing operators of higher derivatives on conventional grids consists of generating spurious saw-tooth solutions while violating smoothness of input data (the even-odd problem, see [19], [18] ). We modify the basic SBP formula by introducing the first-derivative operators defined on the shifted (not symmetric) stencils to avoid the even-odd problem. This modification was announced in [4], [5].

4) To use the proposed schemes within the multiblock approach, we consider discretization of interface transmission conditions at the block boundaries. Application of these more efficient schemes allows us to significantly improve our previous multiblock algorithm [26] for solving seismic modeling problems.

5) We show that the conventional spectral element method (SEM) [10] for second-order equations can be treated as the multiblock finite-difference method whose blocks are the SEM cells with SBP operators. In other words, such SEM algorithm belongs to class of derived schemes within the multiblock framework. The basic observation here is that the differentiation operator of SEM is a special case of SBP operators on nonuniform grids, see [6].

The paper is organized as follows. *Section 2* describes the problem formulation and basic approaches used. In *Section 3*, we consider one-dimensional wave equation. First, we introduce SBP differentiation operators on a uniform grid with symmetric and shifted stencils, and on the SEM cell stencil: Section 3.1. Then, we discretize wave equation with Dirichlet, Robin, and non-reflecting boundary conditions: Section 3.2. In Section 3.3, we introduce the multiblock scheme and correspondent discretization. In Section 3.4, we discuss the accuracy of the proposed schemes. *Section 4* is devoted to approximation of 3D Navier wave equation in curvilinear coordinates. Free surface, Dirichlet, and non-reflecting boundary conditions are considered in Sections 4.1 to 4.3. Section 4.4 contains formulas for the multiblock approach. Stability and accuracy of proposed difference schemes are discussed in Section 4.5. In Section 4.6 we show that the SEM scheme is a combination of SBP and multiblock approaches for the considered equations. *Section 5* contains examples of 1D, 2D (Lamb problem), and 3D test calculations. Conclusions are formulated in *Section 6*.

Note that the SBP difference schemes are generated on the basis of the tensor product of the one-dimensional Green's formula discrete analog, i.e., the canonic computational domains are rectangles or parallelepipeds. For the domains with curvilinear boundaries, the multidimensional analogues of Green's formulas may require. The corresponding theory was founded by V.S. Ryaben'kii in [21] and then developed by him and his colleagues in the framework of the difference potentials method; see [22] and references therein.

## 2  Problem formulation and basic approaches used

Let us consider the Navier wave equation describing propagation of elastic waves in a curvilinear hexahedral domain for the displacement vector $\boldsymbol{u} = (u_1, u_2, u_3)$ with the source function $\boldsymbol{F} = (F_1, F_2, F_3)$:

$$\begin{cases} \rho \dfrac{\partial^2 u_i}{\partial t^2} - \nabla_j \sigma_{ij} = F_i, \quad i = 1, 2, 3 \\ \sigma_{ij} = C_{ijkl} \varepsilon_{kl}, \quad \varepsilon_{kl} = 0.5\ \nabla_k u_l + \nabla_l u_k \end{cases} \quad (1)$$

(hereinafter we use the Einstein summation convention, i.e., sum over the repeated indices in products (including differentiation operators); free indices tell us how many equations have been compressed into one). Here, $\rho$ is the density, $C_{ijkl}$ is the elastic stiffness tensor [14], $i, j, k, l = \{1, 2, 3\}$. We suppose that these coefficients, up to 22 independent ones for general anisotropy, are smooth functions of coordinates $(x_1, x_2, x_3)$. We consider uniform initial conditions at $t = 0$ for simplicity:

$$\boldsymbol{u}\,|_{t=0} = \dfrac{\partial}{\partial t} \boldsymbol{u}\,|_{t=0} = 0. \quad (2)$$

The boundary conditions (BC) on the faces of the computational domain can be either *free surface* BC governed by equations

$$\sigma_{ij} \nu_j = 0 \quad (3)$$

where $(\nu_1, \nu_2, \nu_3)$ is the normal vector, or nonuniform *Dirichlet* BC

$$\boldsymbol{u} = \boldsymbol{a} \quad (4)$$

or, *non-reflecting* BC for open boundaries in the form

$$A_{ij}\frac{\partial u_j}{\partial t}+\sigma_{ij}\nu_j=0 \tag{5}$$

where the matrices $A_{ij}$ are defined hereafter, see (51).

Our aim is to generate finite-difference schemes with high-order spatial accuracy for the initial boundary value problems governed by (1), (2), and conditions from the set (3) - (5). We are also adding the following options that extend the capabilities of the method:

- curvilinearity of computational domain and grid, which is generated by a coordinate transformation of the canonical parametric cube with a rectangular grid;

- multiblock structure of the computational domain, which is provided by explicit treatment of transmission conditions at the block boundaries (continuity of displacements and normal stresses):

$$\boldsymbol{u}=0,\ \left[\sigma_{ij}\nu_j\right]=0 \tag{6}$$

We use the SBP approach for generating high-order finite-difference schemes. To extend application field of the schemes for cases of insufficient smoothness of the equation coefficients and/or point sources for right-hand side (RHS) we modify the basic SPB formula, see (17). This permits us to avoid issues with spurious saw-tooth solutions.

The SBP approach helps to generate symmetric matrices of the approximated spatial operators. This is important property for providing stability of time integration. Below we review the result for hyperbolic systems with the second-order time derivative approximated by the explicit central difference scheme.

In $R^N$ consider a Euclidian space of vector functions $u(t)\in R^N$ with the scalar product $(u,v)=v^{\mathrm{T}}u$ and operators ($N\times N$ matrices) $G$, $A$, etc. (hereinafter «$^{\mathrm{T}}$» denotes the transposition). Define also the scalar product $(u,v)_G=(Gu,v)$ and the corresponding norm $\|u\|_G$ for $G = G^T > 0$. Denote the time derivatives $u_t = \dfrac{du}{dt}$, $u_{tt} = \dfrac{d^2 u}{dt^2}$. We need subsequently the following statement.

**Lemma.** *Suppose for the task*

$$Gu_{tt} + Bu_t + Au = f, \quad u(0) = u_t(0) = 0, \quad u(t), f(t) \in R^N, \ t \geq 0 \qquad (7)$$

*with some operators* $A = A^T > 0$, *diagonal* $G = G^T > 0$, *and diagonal* $B \geq 0$ *satisfying*

$$\lambda_{\min}(Gv, v) \leq (Av, v) \leq \lambda_{\max}(Gv, v), \quad \lambda_{\min} > 0 \qquad (8)$$

*given the explicit difference scheme,* $k = 1, 2, ....,$

$$G\frac{u^{k+1} - 2u^k + u^{k-1}}{\tau^2} + B\frac{u^{k+1} - u^{k-1}}{2\tau} + Au^k = f^k, \quad u^0 = u^1 = 0, \quad u^k, f^k \in R^N \qquad (9)$$

*Then the scheme is stable if*

$$\tau \leq 2 / \sqrt{(1 + \varepsilon)\lambda_{\max}} \qquad (10)$$

*where* $\varepsilon > 0$ *is any number independent of* $\lambda_{\max}$ *and* $\tau$; *the solution is estimated by*

$$\|u^{k+1}\|_G \leq \frac{\text{const}}{\lambda_{\min}} \sqrt{\frac{1 + \varepsilon}{\varepsilon}} \sum_{i=1}^{k} \tau \|f^i\|_G \qquad (11)$$

The lemma follows directly from Theorem 9, Section 6.3 in [24].

The truth of the lemma is sufficient for our analysis of the proposed schemes. We also note that traditionally the energy method is used for the stability study of difference schemes on the basis of SBP approach for second-order equations, see [18] and references therein.

# 3  One-dimensional wave equation

We will generate spatial operators approximating (1) by using the tensor product of one-dimensional operators. The main ideas of the proposed difference schemes are seen from the scalar case. So, let us start with analysis of the one-dimensional wave equation. We consider traditional SBP schemes for uniform grids. We also analyze the difference scheme generated by the spectral elements method, and show that this is a special case of SBP schemes on non-uniform grids inside a single element.

Let us consider the one-dimensional wave equation on the interval $x \in [0,1]$:

$$\rho \frac{\partial^2 u}{\partial t^2} - \frac{\partial}{\partial x} \rho c^2 \frac{\partial}{\partial x} u = f , \qquad (12)$$

where the density $\rho(x)$ and the sound speed $c(x)$ are smooth functions. Adding uniform initial conditions (for simplicity)

$$u \big|_{t=0} = \frac{\partial u}{\partial t} \big|_{t=0} = 0 \qquad (13)$$

and boundary conditions (Dirichlet or Neumann, for example)

$$lu \big|_{x=0} = a_L(t), \quad lu \big|_{x=1} = a_R(t) \qquad (14)$$

we obtain the initial boundary value problem (12)-(14).

## 3.1  *Approximation of the first-derivative operator*

Introduce a grid with $N$ points on the interval $[0,1]$:

$$0 = x_1 < x_2, ..., x_n, ... x_{N-1} < x_N = 1. \qquad (15)$$

Denote by $\bar{u}(t)$ or $\bar{u}$ the trace of a function $u(t,x)$ on the grid (15). Introduce the scalar product

$$(\bar{u}, \bar{v})_H = \bar{u}^T H \bar{v} = \bar{v}^T H \bar{u} \qquad (16)$$

with the weight $H$, where $H > 0$ is a diagonal matrix, and two grid operators $D^+$ and $D^-$, such that

1. $D^+$ and $D^-$ approximate the first derivative $\partial / \partial x$ in all grid points (15)
2. $D^+$, $D^-$, and $H$ satisfy the SBP condition:

$$\left(\overline{u}, D^{+}\overline{v}\right)_H + \left(D^{-}\overline{u}, \overline{v}\right)_H = -\overline{u}\,|_{x=0}\,\overline{v}\,|_{x=0} + \overline{u}\,|_{x=1}\,\overline{v}\,|_{x=1}, \quad \forall \overline{u}, \overline{v} \in R^N. \qquad (17)$$

Eq. (17) is the discrete counterpart of the integration by parts

$$\int_0^1 u \frac{\partial v}{\partial x} dx + \int_0^1 v \frac{\partial u}{\partial x} dx = -uv\big|_{x=0} + uv\big|_{x=1},$$

and a modification of traditional formula with $D \equiv D^+$, $D^- = D^T$, see e.g. [27]. Let us generate such $D^+$, $D^-$, and $H$.

### 3.1.1 Uniform grid, symmetric stencil

Consider the uniform grid (15):

$$x_n = (n-1)h, \quad n = 1, ..., N, \quad h = \frac{1}{N-1}. \qquad (18)$$

In the original SBP approach, the condition (17) is formulated for the same operator $D$ so that $D^+ = D$, $D^- = D^T$. In particular, in [27] such operators $D$ have been derived on the basis of the standard central-difference operators of $p$-th approximation order in the inner grid points of the interval for $p = 2, 4, 6, 8$. The operators have nonsymmetric stencils in $p$ points near the boundary and a smaller, not less than $p/2$, order of accuracy.

These operators work well for approximation of (12)-(14) on smooth input data. However, for insufficiently smooth coefficients, the difference scheme generates saw-tooth waves: the so-called even-odd issue. The propagation velocity of these spurious waves exceeds the physical velocity in the medium and depends on the order of the scheme. In particular, the effect is clearly visible when the acoustic velocity is a stepwise function or when a point source is used in RHS of (12). An example of such solution is given in Section 5.1.

### 3.1.2 Uniform grid, shifted stencils

To avoid the even-odd issue, we use the dual operators $D^+$ and $D^-$ corresponding to the forward and backward differences. The process of their construction is different for internal and near-boundary points. Let 0 be index of the stencil point where the derivative is approximated. Operator $D^+$ for the internal grid points is defined on the asymmetric stencil $[-p/2+1, ..., 0, ..., p/2+1]$ by using the conventional Taylor expansions method with an even approximation order $p$. The matrix $H$ entries are equal to the grid spacing at the interior points. In the near-boundary points, the operators $D^+$, $D^-$, and $H$ are jointly found out to fulfill the conditions (17); the first derivative is approximated with not less than $p/2$ order.

We generated triplets of such $D^+$, $D^-$, and $H$ for $p = 4, 6, 8$. Below is an example for $p = 4$, $N > 8$. The matrix $H$ is diagonal with positive entries:

$$H = \text{diag}\left(\frac{49}{144}, \frac{61}{48}, \frac{41}{48}, \frac{149}{144}, 1, ..., 1, \frac{149}{144}, \frac{41}{48}, \frac{61}{48}, \frac{49}{144}\right) h,$$

where $h$ is the spacing of grid (15). Let us describe the strings $d^+_n$ and $d^-_n$, $n = 1, ..., N$ of matrices $D^+$ and $D^-$, respectively. The first $p$ stings of $D^+$ are the following:

$d^+_1 = [-59/42, 12/7, -3/14, -2/21, 0, 0, 0, 0, ..., 0]/h$

$d^+_2 = [-103/183, 15/122, 31/61, -49/366, 4/61, 0, 0, 0, ..., 0]/h$

$d^+_3 = [59/246, -38/41, -21/82, 176/123, -24/41, 4/41, 0, 0, ..., 0]/h$

$d^+_4 = [-5/447, 15/298, -51/149, -665/894, 216/149, -72/149, 12/149, 0, ..., 0]/h$

The first $p$ stings of $D^-$:

$d^-_1 = [-451/294, 103/49, -59/98, 5/147, 0, 0, 0, 0, ..., 0]/h$

$d^-_2 = [-28/61, -15/122, 38/61, -5/122, 0, 0, 0, 0, ..., 0]/h$

$d^-_3 = [7/82, -31/41, 21/82, 17/41, 0, 0, 0, 0, ..., 0]/h$

$d^-_4 = [14/447, 49/298, -176/149, 665/894, 36/149, 0, 0, 0, ..., 0]/h$

At the internal points:

$$\begin{cases} d^+_n(n-1:n+3) = [-1/4, -5/6, 3/2, -1/2, 1/12]/h \\ d^+_n(1:n-2) = d^+_n(n+4:N) = 0 \end{cases}, n = 5, ..., N-4$$

$$\begin{cases} d^-_n(n-3:n+1) = [-1/12, 1/2, -3/2, 5/6, 1/4]/h \\ d^-_n(1:n-4) = d^-_n(n+2:N) = 0 \end{cases}, n = 5, ..., N-4 \quad (19)$$

The last $p$ strings of $D^+$ and $D^-$ are calculated according to the rule:

$d^+_i(j) = -d^-_{N-i+1}(N-j+1), \forall i, j = 1, N$;

for instance

$d^+_N = -d^-_1(N:-1:1) = [0, ..., 0, 0, 0, 0, -5/147, 59/98, -103/49, 451/294]/h$.

The advantage of $D^+$ and $D^-$ over the conventional operators $D^+ = D$, $D^- = D^T$ is shown in Section 5.1.

*3.1.3 Nonuniform grid, SEM cell stencil*

An interesting example of SBP operators gives the SEM, see [6] where this observation is discussed for the discontinuous Galerkin SEM approach. Let us consider the interval $[-1;1]$ as a single SEM cell. The GLL grid is used in (15), i.e., the nodes $x_n$, $n = 1,...,N$ are the roots of the polynomial $(1-x^2)P'_{N+1}(x)$ where $P_{N+1}$ is the Legendre polynomial of $(N+1)$-th degree. The matrix $H$ is defined by

$$H = \text{diag}(\omega_n),$$
$$\omega_n = \frac{2}{(N+1)(N+2)P_{N+1}^2(x_n)}, \; n = 2,...,N-1; \; \omega_n = \frac{2}{(N+1)(N+2)}, \; n = 1, N,$$

i.e., it provides the following Gaussian quadrature formula:

$$\int_{-1}^{1} u \, dx \approx \sum_{n=1}^{N} u(x_n) \omega_n \qquad (20)$$

This formula is exact for all polynomials up to the degree $2N-3$. This important property is a link from the SEM to the SBP approach, see (22), (23). The operator $D$ of the first derivative is the full $N \times N$ matrix defined by $D_{nm} = l'_m(x_n)$ where

$$l_m(x) = \prod_{j=1, j \neq m}^{N} \frac{x - x_m}{x_j - x_m}$$

are the interpolation Lagrange polynomials, and $l_m(x_n) = \delta_{mn}$, $\delta_{mn}$ is the Kronecker symbol.

One can easily see that the operators $H$, $D^+ = D$, $D^- = D^T$ satisfy conditions 1, 2 from (17). Indeed, item 1 is valid by construction. The condition 2

$$\langle \overline{u}, D\overline{v} \rangle_H + \langle D^T \overline{u}, \overline{v} \rangle_H = -\overline{u}|_{x=0} \overline{v}|_{x=0} + \overline{u}|_{x=1} \overline{v}|_{x=1}, \; \forall \overline{u}, \overline{v} \in R^N \qquad (21)$$

is proved as follows. In the identity

$$\int_{-1}^{1} l'_m(x) l_k(x) dx + \int_{-1}^{1} l_m(x) l'_k(x) dx = l_m(1) l_k(1) - l_m(-1) l_k(-1) = \delta_{mN} \delta_{kN} - \delta_{m1} \delta_{k1}$$

we can use the equations

$$\int_{-1}^{1} l'_m(x) l_k(x) dx = \sum_{n=1}^{N} D_{nm} l_{kn} \omega_n = D_{km} \omega_k = (HD)_{km} \qquad (22)$$

$$\int_{-1}^{1} l_m(x) l'_k(x) dx = \sum_{n=1}^{N} l_{mn} D_{nk} \omega_n = D_{mk} \omega_m = (D^T H)_{mk} \qquad (23)$$

due to exact quadratures (20) for polynomials up to the ($2N-3$)-th degree. Therefore,

$$(D^T H)_{mk} + (HD)_{km} = \delta_{mN} \delta_{kN} - \delta_{m1} \delta_{k1},$$

i.e.,

$$D^T H + HD = Q,$$

where $Q = \mathrm{diag}(-1, 0, ..., 0, 1)$; this proves (21).

## 3.2 Discretization of the problem (12)-(14)

To satisfy conditions of the lemma, we generate a symmetric matrix of the spatial operator during discretization.

Define the positive diagonal matrices $P = \mathrm{diag}(\{\rho(x_n)\})$ and $C = \mathrm{diag}(\{c(x_n)\})$ with values of density and velocity on the grid points (15). Denote the time derivatives of a vector function $\bar{u}$ given on this grid by $\bar{u}_t$, $\bar{u}_{tt}$. Suppose we have the operators $D^+$, $D^-$, $H$ satisfying the conditions 1, 2, see (17).

Let us take arbitrary $\bar{w} \in R^N$ and substitute $\bar{u} = PC^2 D^+ \bar{w}$ for (17); we have

$$\begin{aligned} &\langle PC^2 D^+ \bar{w}, D^+ \bar{v} \rangle_H + \langle D^- PC^2 D^+ \bar{w}, \bar{v} \rangle_H = \\ &= -(PC^2 D^+ \bar{w})|_{x=0}\, \bar{v}|_{x=0} + (PC^2 D^+ \bar{w})|_{x=1}\, \bar{v}|_{x=1}, \quad \forall \bar{w}, \bar{v} \in R^N \end{aligned} \qquad (24)$$

It follows from (24) and (16) that

$$\bar{v}^T (D^+)^T H P C^2 D^+ \bar{w} + \bar{v}^T H D^- P C^2 D^+ \bar{w} = \bar{v}^T Q P C^2 D^+ \bar{w}, \quad \forall \bar{w}, \bar{v} \in R^N.$$

Hence, the operator $HD^- PC^2 D^+ - QPC^2 D^+$ has a symmetric matrix. Evidently, a constant function is its kernel. Hence, due to positivity of $P$ and $C$, the matrix is negatively defined in the subspace $R^N \setminus \{\mathrm{const}\}$. Therefore, in accordance with the lemma, the time integration (9) of the system

$$HP\bar{u}_{tt} - (HD^- PC^2 D^+ - QPC^2 D^+)\bar{u} = H\bar{f} \qquad (25)$$

with this spatial operator will be stable. Denote

$$\tilde{L} \equiv HD^-PC^2D^+ - QPC^2D^+ = -(D^+)^T HPC^2D^+.$$

Since the operators $D^-$, $D^+$ approximate the first derivative the semidiscrete system

$$P\bar{u}_{tt} - H^{-1}\tilde{L}\bar{u} = \bar{f} \qquad (26)$$

approximates the wave equation (12) in the internal points $\{x_2,...,x_{N-1}\}$ of the interval $[0,1]$.

Now let us consider approximation of boundary conditions.

### 3.2.1 Dirichlet conditions

Take $lu \equiv u$ in (14). Represent a solution in the form of $\bar{u} = \bar{u}^0 + \bar{d}$ where

$$\bar{u}^0 = \{0, u_2,...,u_{N-1}, 0\}^T$$

$$\bar{d} = \{a_L, 0,...,0, a_R\}^T,$$

and consider (26) for $\bar{u}^0$ only. To formalize this we introduce the operator $R = \text{diag}(0,1,...,1,0)$ projecting any $\bar{u}$ into $\bar{u}^0$. Due to $RQ = 0$ we have

$$\begin{aligned}R\tilde{L}\bar{u} &= RHD^-PC^2D^+\bar{u} = RHD^-PC^2D^+\bar{u}^0 + RHD^-PC^2D^+\bar{d}\\ &\approx RHP\bar{u}^0_{tt} - RH\bar{f} + RHD^-PC^2D^+\bar{d}\end{aligned} \qquad (27)$$

Hereafter, the sign "$\approx$" is used while substituting *a solution* of (12)-(14) for the considered formulas, i.e., "approximation on the solution" is considered.
Therefore, we formulate the following problem for Dirichlet conditions in (14):

$$RHPR\bar{u}^0_{tt} - RHD^-PC^2D^+R\bar{u}^0 = RHR\bar{f} - RHD^-PC^2D^+\bar{d}, \quad u(0) = u_t(0) = 0 \ . \quad (28)$$

It has a symmetric matrix in the subspace $\bar{u}^0 \in R^{N-2}$ of unknown functions with zeros at the interval ends. Therefore, in accordance with the lemma we obtain the following theorem.

**Theorem 1**. *The explicit scheme (9) of the time integrating the Dirichlet problem (28) with the diagonal matrix $G = RHP$ and the symmetric negatively defined matrix $A = RHD^-PC^2D^+R$ is stable under the condition (10).*

### 3.2.2 Robin conditions

Take $lu \equiv \dfrac{\partial u}{\partial x} + b_L(t)u$ at $x=0$ and $lu \equiv \dfrac{\partial u}{\partial x} + b_R(t)u$ at $x=1$ in (14). Here

$$b_L(t) \leq 0, \quad b_R(t) \geq 0 \qquad (29)$$

are given smooth functions. Denote them by a single function $b = b_L, b_R$ and introduce the operator $lu \equiv \dfrac{\partial u}{\partial x} + b(t)u$ for brevity.

We have

$$\tilde{L}\bar{u} = HD^-PC^2D^+\bar{u} - QPC^2D^+\bar{u} \approx HD^-PC^2D^+\bar{u} - QPC^2(\bar{d} - b\bar{u})$$
$$\approx HP\bar{u}_{tt} - H\bar{f} - QPC^2\bar{d} + bQPC^2\bar{u}$$

Therefore, we formulate the following problem with the symmetric matrix in the space $\bar{u} \in R^N$ of unknown functions for the Robin condition in (14):

$$HP\bar{u}_{tt} - (\tilde{L} - bQPC^2)\bar{u} = H\bar{f} + QPC^2\bar{d} . \qquad (30)$$

As we can see, if at least one of the functions $b(t)$ is not equaled to zero in the condition (29) then $\tilde{L} - bQPC^2$ is the sign-definite operator. The case $b(t) \equiv 0$ corresponds to the Neumann conditions, so, as usual, we need to exclude the constant function from solutions space by an additional condition. In accordance with the lemma we obtain the following theorem.

**Theorem 2**. *The explicit scheme (9) of the time integrating the Robin problem (30) with the diagonal matrix $G = HP$ and the symmetric negatively defined matrix $A = \tilde{L} - bQPC^2$ (in the subspace $R^N \setminus \{\text{const}\}$ if $b(t) \equiv 0$) is stable under the condition (10).*

### 3.2.3 Non-reflecting conditions

Take $lu \equiv \dfrac{\partial u}{\partial t} + b\dfrac{\partial u}{\partial x}$ in (14) with $b = b_L, b_R$ similarly to the Robin problem, and with $a_L = a_R = 0$. We have

$$\tilde{L}\bar{u} = HD^-PC^2D^+\bar{u} - QPC^2D^+\bar{u} \approx HD^-PC^2D^+\bar{u} + b^{-1}QPC^2\dfrac{\partial \bar{u}}{\partial t}$$
$$\approx HP\bar{u}_{tt} - H\bar{f} + b^{-1}QPC^2\dfrac{\partial \bar{u}}{\partial t}$$

Therefore, we formulate the following problem with the symmetric matrix in the space $\bar{u} \in R^N$ of unknown functions for the non-reflecting condition in (14):

$$H\mathrm{P}\bar{u}_{tt} - \tilde{L}\bar{u} + b^{-1}Q\mathrm{P}C^2 \frac{\partial \bar{u}}{\partial t} = H\bar{f}. \qquad (31)$$

Imposing the inequalities

$$b_L(t) \leq 0, \quad b_R(t) \geq 0, \quad b = \{b_L, b_R\} \qquad (32)$$

we satisfy the conditions of the lemma. Therefore, we obtain the following theorem.

**Theorem 3**. *The explicit scheme (9) of the time integrating the problem (31) with the diagonal matrix $G = H\mathrm{P}$, matrix $B = b^{-1}Q\mathrm{P}C^2$, and the symmetric negatively defined matrix $A = \tilde{L}$ in the subspace $R^N \setminus \{const\}$ is stable under the conditions (10) and (32).*

**Remark.** The formulas derived above, e.g., (30), for the Robin problem, coincide with those that can be obtained if one follows the SBP-SAT approach [16]. Indeed, the symmetry and negative definiteness of the matrix $\tilde{L}$ allows us to define the integral $E(t) \equiv 0.5(Gu_t, u_t) + 0.5(Au, u)$ of the system (7), see the lemma. In turn, this is the integral that is used in SBP-SAT approach to derive governing equations.

## 3.3 Multiblock approach

The multiblock approach is useful if the density $\rho(x)$ and sound velocity $c(x)$ have strong discontinuities at some points. Let us consider (12) with smooth coefficients to the left and to the right from $x = 0$. As usual, we impose the following transmission conditions for the case of discontinuous coefficients at $x = 0$:

$$\begin{aligned} u_{0-}(t,0) &= u_{0+}(t,0), \\ \rho_{0-} c_{0-}^2 \left.\frac{\partial u_{0-}}{\partial x}\right|_{x=0} &= \rho_{0+} c_{0+}^2 \left.\frac{\partial u_{0+}}{\partial x}\right|_{x=0}. \end{aligned} \qquad (33)$$

We will use the symbols $(.)|_{0+}$ and $(.)|_{0-}$ to indicate that the value at $x=0$ is calculated by difference formulas from the right and from the left, respectively. We have

$$\tilde{L}\bar{u}\,|_{0+} = HD^- PC^2 D^+ \bar{u}\,|_{0+} - QPC^2 D^+ \bar{u}\,|_{0+} =$$
$$\approx HP\bar{u}_{tt}\,|_{0+} - H\overline{f}\,|_{0+} + PC^2 D^+ \bar{u}\,|_{0+}\,,$$

and similarly:

$$\tilde{L}\bar{u}\,|_{0-} \approx HP\bar{u}_{tt}\,|_{0-} - H\overline{f}\,|_{0-} - PC^2 D^+ \bar{u}\,|_{0-}.$$

Let us sum these equations in the common point $x=0$ and take into account (33). We have

$$0 = HP\bar{u}_{tt}\,|_{0+} + HP\bar{u}_{tt}\,|_{0-} - \tilde{L}\bar{u}\,|_{0+} - \tilde{L}\bar{u}\,|_{0-} +$$
$$+ PC^2 D^+ \bar{u}\,|_{0+} - PC^2 D^+ \bar{u}\,|_{0-} - H\overline{f}\,|_{0+} - H\overline{f}\,|_{0-} \approx$$
$$\approx HP\bar{u}_{tt}\,|_{0+} + HP\bar{u}_{tt}\,|_{0-} - \tilde{L}\bar{u}\,|_{0+} - \tilde{L}\bar{u}\,|_{0-} - H\overline{f}\,|_{0+} - H\overline{f}\,|_{0-}\,.$$

Therefore, we formulate the following equation for $x=0$:

$$HP\bar{u}_{tt}\,|_{0+} + HP\bar{u}_{tt}\,|_{0-} - \tilde{L}\bar{u}\,|_{0+} - \tilde{L}\bar{u}\,|_{0-} = H\overline{f}\,|_{0+} + H\overline{f}\,|_{0-}.$$

As the equations in the grid points with $x<0$ and with $x>0$ have the form

$$HP\bar{u}_{tt} - \tilde{L}\bar{u} = H\overline{f}\,,$$

we see that the symmetry of the matrix for spatial operator is preserved; hence, the time integration by the scheme (9) is stable.

**Remark.** Although the notation for the operators in the above formulas is same in cases $(.)|_{0-}$ and $(.)|_{0+}$, the grids and operators themselves for $x<0$ and $x>0$ can be different, of course.

### 3.4 On the accuracy of the proposed difference schemes

Analysis of the time-step upper limit in (10), as well as use of (11) to evaluate the convergence rate, requires estimates of the maximum and minimum eigenvalues in (8). Since the operator $D^- PC^2 D^+$ approximates $\frac{d}{dx}(\rho c^2 \frac{d}{dx})$, the minimum eigenvalue is bounded from below by a positive constant independent of the number of grid points $N$ and the grid type (uniformly distributed or GLL nodes). Therefore, the numerical solution accuracy will be determined by the approximation error of the difference schemes on solutions to the original differential problem. Although the operators $D^+$, $D^-$ have the formal order $p/2$ at the grid points near the boundary (in the uniform grid case), it is possible to optimize their coefficients to improve the accuracy for a certain range of wave numbers in the spirit of [8]. Therefore, it is very important to evaluate numerically the real accuracy and so-called numerical order of proposed difference schemes using, for example, artificial analytical solutions (recall that in this method one solves the problems with the RHS obtained after substitution of a given function for differential operator of the problem). Another important characteristic of high-order schemes is the solution behavior when the medium parameters and/or the RHS are not sufficiently smooth; for example, piecewise smooth density and sound velocity, point sources, etc. Analysis of these properties without corresponding numerical experiments cannot be complete.

Obtaining estimates of the theoretical maximum and minimum eigenvalues of problems (28), (30), and (31) is a separate task and not a goal of this paper. Evidently, the maximum eigenvalue that determines the maximal time integration step depends on the grid type; its magnitude is $O(N^2)$ and $O(N^4)$ for uniformly distributed and GLL nodes, respectively. We only note that with explicit formulas of operators $D^+$, $D^-$, $H$ these estimates are easily obtained numerically with high accuracy during discretization of relevant tasks. For example, numerical limits of $\lambda_{max}/N^2$ as $N \to \infty$ in (8) for the Neumann problem with $G = H$ and $A = HD^-D^+ - QD^+$ are close to 7.1, 4.8, and 5.3 for $p = 4$, $p = 6$, and $p = 8$, respectively (periodical case gives 7.1, 4.8, and 4.4).

Anyway, as usual, the stability is the main problem when building difference schemes; the subsequent experimental investigation of the accuracy of a given difference scheme is much easier.

# 4  Approximation of the Navier wave equation in curvilinear coordinates for heterogeneous anisotropic medium

Application of the approaches described in the previous sections to the elastic multidimensional equations is similar in conceptual level and is made directly. Therefore, we give only the basic equations of the proposed finite-difference spatially highly accurate algorithm of integrating Navier wave equation in curvilinear coordinates. Besides the scheme inside the computational domain, we consider approximation of free surface boundary conditions and propose artificial boundary conditions at the open boundaries to simulate wave propagation without significant reflections.

Consider the system (1). After using a smooth, general-type coordinate transformation

$$x(\xi), \quad x = (x_1, x_2, x_3), \quad \xi = (\xi_1, \xi_2, \xi_3)$$

it can be written in the conservation law

$$\begin{cases} J\rho \dfrac{\partial^2 u_i}{\partial t^2} - \partial_k \ JT_{kj}\sigma_{ij} = F_i(\xi)g(t), \quad i = 1,2,3 \\ \sigma_{ij} = C_{ijkl}\varepsilon_{kl}, \quad \varepsilon_{kl} = 0.5 \ T_{ik}\partial_i u_l + T_{jl}\partial_j u_k \end{cases} \quad (34)$$

where $\partial_k = \dfrac{\partial}{\partial \xi_k}$, $T_{kj} \equiv \dfrac{\partial \xi_k}{\partial x_j}$, $J = \det \dfrac{\partial(x_1, x_2, x_3)}{\partial(\xi_1, \xi_2, \xi_3)}$, $i, j, k, l = 1,2,3$.

For the isotropic 2D case the form (34) is used in [1].

We introduce the computational domain $\xi_{i,\min} \leq \xi_i \leq \xi_{i,\max}$, $i = 1,2,3$, in parametric coordinates, and generate a rectangular grid with $N_i$ nodes in each direction

$$\xi_{i \ n_i}, \quad n_i = 1,...,N_i, \ i = 1,2,3 \quad (35)$$

that can be either with uniform nodes

$$\xi_{i\ n_i} = \xi_{i,\min} + (n_i - 1)h_i, \quad h_i = (\xi_{i,\max} - \xi_{i,\min})/(N_i - 1) = \text{const}$$

or with GLL nodes. Let us consider the following objects on (35) – arrays of length 1, 3, and 9:

$u_i$, the scalar grid functions

$\{u_i\} = (u_1, u_2, u_3)$, the vector grid functions

$\{u_{ij}\}$, the 3x3 tensor grid functions.

Denote by $J$, $T_{kj}$, $C_{ijkl}$ the diagonal $N_1 N_2 N_3 \times N_1 N_2 N_3$ matrix operators with coefficients from (34) for the scalar grid functions; 1, 9, and 81 operators, respectively.

We define the following one-dimensional operators for each direction $\xi_i$ on grid (35) with $N_i \times N_i$ matrices: the identity operators $I^i$; SBP operators $D^{i,+}$, $D^{i,-}$; the weight operators $H^i$ with diagonal matrices; and operators $Q^i = \text{diag}(-1, 0, ..., 0, 1)$. Using tensor products, we construct the desired operators for scalar functions on a 3D grid (35):

$D_i^+$ and $D_i^-$ for approximating the first derivatives $\partial_i$, $i = 1, 2, 3$, along direction $\xi_i$ in (34),

$$D_1^\alpha = D^{1,\alpha} \otimes I^2 \otimes I^3, \quad D_2^\alpha = I^1 \otimes D^{2,\alpha} \otimes I^3, \quad D_3^\alpha = I^1 \otimes I^2 \otimes D^{3,\alpha}, \quad \alpha = \pm,$$

$H$ with the diagonal matrix

$$H = H^1 \otimes H^2 \otimes H^3,$$

and $Q_i$ with diagonal matrices

$$Q_1 = Q^1 \otimes H^2 \otimes H^3, \quad Q_2 = H^1 \otimes Q^2 \otimes H^3, \quad Q_3 = H^1 \otimes H^2 \otimes Q^3$$

Let us introduce the scalar product for scalar functions

$$(u_i, u_j)_H \equiv (Hu_i, u_j) \equiv (u_j)^T H u_i.$$

For vector and tensor functions the scalar products are defined by summing of three components

$$(\{u_i\},\{v_i\})_H \equiv (u_i, v_i)_H, \quad (\{u_i\},\{v_i\}) \equiv (u_i, v_i)$$

and nine components, respectively:

$$(\{u_{ij}\},\{v_{ij}\})_H \equiv (u_{ij}, v_{ij})_H, \quad (\{u_{ij}\},\{v_{ij}\}) \equiv (u_{ij}, v_{ij}).$$

The multidimensional counterpart of SBP formula(17) for the scalar functions reads:

$$u_i, D_k^+ u_j \big|_H + D_k^- u_i, u_j \big|_H = H^{-1} Q_k u_i, u_j \big|_H, \quad i, j, k = 1, 2, 3. \tag{36}$$

Summation of nine combinations of (36) gives the following multidimensional SBP formula for vector and tensor functions:

$$\{u_{ki}\}, \{D_k^+ v_i\} \big|_H + \{D_k^- u_{ki}\}, \{v_i\} \big|_H = \{H^{-1} Q_k u_{ki}\}, \{v_i\} \big|_H. \tag{37}$$

Denote by $U^N$ the space of vector functions $\{u_i\}$ on the grid (35) excluding constants $\{\text{const}_i\}$, the kernel of $D_i^+$, $D_i^-$. Dimension of $U^N$ is $3N_1 N_2 N_3 - 3$.

By analogy with (24) we take a grid vector function $\{w_i\} \in U^N$ and firstly calculate the strain tensor function $\{\varepsilon_{kl}\} = \{0.5 \ T_{ik} D_i^+ w_l + T_{jl} D_j^+ w_k\}$. Afterwards, we calculate the stress tensor function $\{\sigma_{ij}\} = \{C_{ijkl} \varepsilon_{kl}\} = \{C_{ijkl} 0.5 \ T_{i'k} D_{i'}^+ w_l + T_{j'l} D_{j'}^+ w_k\} = \{C_{ijkl} T_{i'k} D_{i'}^+ w_l\}$ taking into account $C_{ijkl} = C_{ijlk}$, and finally obtain the tensor function

$$\{u_{ki}\} = \{JT_{kj} \sigma_{ij}\} = \{JT_{kj} C_{ijk'l} T_{i'k'} D_{i'}^+ w_l\}.$$

Substitute $\{u_{ki}\}$ for (37):

$$\{JT_{kj} C_{ijk'l} T_{i'k'} D_{i'}^+ w_l\}, \{D_k^+ v_i\} \big|_H + \{D_k^- JT_{kj} \sigma_{ij}\}, \{v_i\} \big|_H = \{H^{-1} Q_k JT_{kj} \sigma_{ij}\}, \{v_i\} \big|_H. \tag{38}$$

Let us consider the first term in (38). It consists of the sum of nine scalar products of scalar grid functions. Transposing $D_k^+$, transferring from left to right the coefficients of diagonal operators $J$, $H$, $T_{kj}$, $C_{ijk'l}$ in these scalar products, rearranging them, and using symmetry $C_{lk'ji} = C_{jilk'} = C_{ijk'l}$ we obtain the following chain of equalities

$$\{JT_{kj}C_{ijk'l}T_{i'k}D_{i'}^+w_l\}, \{D_k^+v_i\}_H = \{(D_k^+)^T HJT_{kj}C_{ijk'l}T_{i'k}D_{i'}^+w_l\}, \{v_i\}$$

$$= \{w_l\}, \{(D_{i'}^+)^T T_{i'k}C_{ijk'l}T_{kj}JHD_k^+v_i\}$$

$$\equiv \{w_i\}, \{(D_k^+)^T T_{kj}C_{lk'ji}T_{i'k}JHD_{i'}^+v_l\} = \{w_i\}, \{(D_k^+)^T HJT_{kj}C_{ijk'l}T_{i'k}D_{i'}^+v_l\}$$

(the sign of identity is used when we simply rename indices). Thus the operator in equation $\{u_i\} := \{(D_k^+)^T HJT_{kj}C_{ijk'l}T_{i'k}D_{i'}^+w_l\}$ has a symmetric matrix. Due to positiveness of operators $J$, $H$, and the quadratic form $(\{C_{ijkl}u_{kl}\}, \{u_{ij}\})$, see [14], this matrix is positively defined for grid vector functions belonging $U^N$.

Hence, we conclude from (38) that the operator in equation

$$\{u_i\} := \{(HD_k^- - Q_k)JT_{kj}C_{ijk'l}T_{i'k}D_{i'}^+w_l\} \tag{39}$$

has a symmetric negatively defined matrix on the space $U^N$.

Therefore, we discretize (34) on the grid (35) by the following ODE system with respect to time:

$$\begin{cases} JP\dfrac{d^2\{u_i\}}{dt^2} - \{(D_k^-JT_{kj} - H^{-1}Q_kJT_{kj})\sigma_{ij}\} = \{F_i\}g(t) \\ \{\sigma_{ij}\} = \{C_{ijj'i'}\varepsilon_{j'i'}\}, \quad \{\varepsilon_{j'i'}\} = \{0.5\ T_{ij'}D_i^+u_{i'} + T_{ji'}D_j^+u_{j'}\} \end{cases} \tag{40}$$

These equations approximate (34) in the inner grid points, i.e., where $Q_k$ vanishes. To investigate the approximation at the boundary points, we should consider the boundary conditions; this will be done below.

The system (40) is simply transformed to a discrete Navier wave equation counterpart:

$$JP\dfrac{d^2\{u_i\}}{dt^2} - \{(D_k^- - H^{-1}Q_k)JT_{kj}C_{ijj'i'}T_{k'j'}D_k^+u_{i'}\} = \{F_i\}g(t) \tag{41}$$

that is written in the form (due to (38)):

$$HJP\dfrac{d^2\{u_i\}}{dt^2} + \{(D_k^+)^T HJT_{kj}C_{ijj'i'}T_{k'j'}D_k^+u_{i'}\} = H\{F_i\}g(t)\ . \tag{42}$$

For brevity, we introduce the notation of the spatial operator in (42), (41):

$$E_{ii'} \equiv -(D_k^+)^T HJT_{kj} C_{ijj'i'} T_{k'j'} D_{k'}^+ = (HD_k^- - Q_k) JT_{kj} C_{ijj'i'} T_{k'j'} D_{k'}^+ \ . \tag{43}$$

Notice that sequential use of $D_i^+$, $D_i^-$ over the whole computational domain in (40) permits us to a) optimize the theoretical number of operations per grid point; b) reduce the amount of memory in the implementation of Hooke's law, and c) build an efficient parallel algorithm for multicore computing systems. For instance, such difference scheme has $O(p)$ of float operations per grid point while computing the mixed derivatives (instead of $O(p^2)$ in case of direct approximation of mixed derivatives by a second-order operator) where $p+1$ is the stencil length.

Time discretization of (42) is made by the standard finite-difference explicit scheme, see (9).

Now consider task formulations with different boundary conditions.

## *4.1 Free surface boundary*

Suppose that the free surface corresponds to the maximal (or minimal) value of $\xi_k$ for any $k = 1, 2, 3$ in parametric coordinates. The normal to it is defined by

$$\{v_{kj}\} = T_{(k)}^{-1} \{T_{kj}\} \tag{44}$$

where

$$T_{(k)} = \begin{bmatrix} -\sqrt{T_{k1}^2 + T_{k2}^2 + T_{k3}^2}, & \xi_k = \xi_{k,\min} \\ \sqrt{T_{k1}^2 + T_{k2}^2 + T_{k3}^2}, & \xi_k = \xi_{k,\max} \end{bmatrix}$$

(index in parentheses indicates that this is a free index without summation).

Thus boundary condition (3) reads

$$0 = \{v_{kj} \sigma_{ij}\} = T_{(k)}^{-1} \{T_{kj} \sigma_{ij}\}, \quad k = 1, 2, 3 \ . \tag{45}$$

Denote by $\{\tilde{u}_i\}$ the trace on grid (35) of a solution to (1), (2) with conditions (3) at all six boundaries. Substituting $\{\tilde{u}_i\}$ for the spatial term of (41) we have

$$\{(-D_k^- + H^{-1}Q_k)JT_{kj}C_{ijj'i'}T_{k'j'}D_{k'}^+\tilde{u}_{i'}\} \approx -JP\frac{d^2\{\tilde{u}_i\}}{dt^2} + \{F_i\}g(t) + \{H^{-1}Q_kJT_{kj}C_{ijj'i'}T_{k'j'}D_{k'}^+\tilde{u}_{i'}\}$$
$$\approx -JP\frac{d^2\{\tilde{u}_i\}}{dt^2} + \{F_i\}g(t)$$

(the first and second signs «$\approx$» correspond to the approximation of equations (1) in all grid points and approximation of conditions (45), respectively). Therefore, using (42) we obtain the problem

$$HJP\frac{d^2\{u_i\}}{dt^2} + \{(D_k^+)^T HJT_{kj}C_{ijj'i'}T_{k'j'}D_{k'}^+ u_{i'}\} = H\{F_i\}g(t) \qquad (46)$$

that approximates (40), (2), and (3) for all boundaries with the symmetric positively defined matrix of the spatial operator on $U^N$.

## 4.2 Dirichlet conditions

Denote by $\{\tilde{u}_i\}$ the trace on grid (35) of a solution to (1), (2) with conditions (4) at all six boundaries. Substituting $\{\tilde{u}_i\}$ for the spatial term of (41) we have

$$\{(-D_k^- + H^{-1}Q_k)JT_{kj}C_{ijj'i'}T_{k'j'}D_{k'}^+\tilde{u}_{i'}\} \approx$$
$$\approx -JP\frac{d^2\{\tilde{u}_i\}}{dt^2} + \{F_i\}g(t) + \{H^{-1}Q_kJT_{kj}C_{ijj'i'}T_{k'j'}D_{k'}^+ a_{i'}\} \qquad (47)$$

Introduce the subspace of all functions in $U^N$ having zero values if the grid index is outside the set $(2,...,N_1-1) \otimes (2,...,N_2-1) \otimes (2,...,N_3-1)$. Denote by $R$ the projection operator of a grid function onto this subspace. Regarding that $RQ_k = 0$, we obtain the following task from (47) and (38):

$$RHJPR\frac{d^2\{u_i\}}{dt^2} + \{R(D_k^+)^T HJT_{kj}C_{ijj'i'}T_{k'j'}D_{k'}^+ Ru_{i'}\} =$$
$$= RH\{F_i\}g(t) - \{RQ_kJT_{kj}C_{ijj'i'}T_{k'j'}D_{k'}^+ a_{i'}\} \qquad (48)$$

It approximates (40), (2) with conditions (4) and has a symmetric positively defined matrix of the spatial operator when solution is sought among functions with the grid indices from the set $(2,...,N_1-1) \otimes (2,...,N_2-1) \otimes (2,...,N_3-1)$.

### 4.3 Open boundaries

The condition (5) imitating open boundaries after the coordinate transformation reads

$$\{A_{k,ij} \frac{du_j}{dt}\} + \{T_{kj}\sigma_{ij}\} = 0, \quad k=1,2,3, \tag{49}$$

where $A_{k,ij} = T_{(k)}A_{ij}$ are some matrices corresponding to six faces of the parametric parallelepiped (two faces for each $k$).

Denote by $\{\tilde{u}_i\}$ the trace on grid (35) of a solution to (1), (2) with conditions (49) at all six boundaries. Substituting $\{\tilde{u}_i\}$ for the spatial term of (41) we have

$$\{(-D_k^- + H^{-1}Q_k)JT_{kj}C_{ijj'i'}T_{k'j'}D_{k'}^+\tilde{u}_{i'}\} \approx -JP\frac{d^2\{\tilde{u}_i\}}{dt^2} + \{F_i\}g(t) - \{H^{-1}Q_k JA_{k,ij}\frac{d\tilde{u}_j}{dt}\}.$$

We obtain the problem

$$HJP\frac{d^2\{u_i\}}{dt^2} + \{Q_k JA_{k,ij}\frac{du_j}{dt}\} + \{(D_k^+)^T HJT_{kj}C_{ijj'i'}T_{k'j'}D_{k'}^+u_{i'}\} = H\{F_i\}g(t) \tag{50}$$

with a symmetric positively defined matrix of the spatial operator on $U^N$ and with non-reflecting conditions at all boundaries.

To satisfy the condition $B \geq 0$ of the lemma we need non-negatively defined matrices of the operator $Q_k JA_{k,ij}$ in (50). Let us consider the following example. We introduce the matrix of normal components

$$N_k = \begin{pmatrix} v_{k1} & 0 & 0 & 0 & v_{k3} & v_{k2} \\ 0 & v_{k2} & 0 & v_{k3} & 0 & v_{k1} \\ 0 & 0 & v_{k3} & v_{k2} & v_{k1} & 0 \end{pmatrix}$$

and 6×6 matrix $\tilde{C}$ of Hooke's law in the matrix (Voigt) notation for stress vector $(\sigma_{11},\sigma_{22},\sigma_{33},\sigma_{23},\sigma_{13},\sigma_{12})^T$ and strain vector $(\varepsilon_{11},\varepsilon_{22},\varepsilon_{33},2\varepsilon_{23},2\varepsilon_{13},2\varepsilon_{12})^T$. Define $A_{k,ij}$ as follows

$$\{A_{k,ij}\} = T_{(k)}\rho(N_k \tilde{C} N_k^T)^{\frac{1}{2}}. \qquad (51)$$

This matrix provides positiveness of $Q_k J A_{k,ij}$. Moreover, plane waves moving in $\{v_{k1}, v_{k2}, v_{k3}\}$ direction in the homogeneous medium defined by $\rho$ and $\tilde{C}$ satisfy (49). For these waves (and only), the boundary is transparent; this corresponds to partial passage of waves without reflection from a boundary. A partial case of matrix (51) for isotropic media was derived in [20].

Let us mention also another local non-reflecting condition obtained in [25] for general second-order hyperbolic systems. Its formulas contain spatial derivatives to take into account heterogeneity of medium parameters. Generally speaking, this condition does not have the form (49) when applied to equations of elasticity.

## *4.4 Multiblock approach*

Consider special case of the boundary between two grid blocks, Eq. (41) is satisfied in each of them. Let the boundary $\Gamma$ between blocks correspond to index $l$ with any value from 1 to 3. Normal to it $\{v_{lj}\} = T_{(l)}^{-1}\{T_{lj}\}$. Boundary conditions of continuity solutions and normal stresses are given at $\Gamma$:

$$\{u_i\}|_{\Gamma^+} = \{u_i\}|_{\Gamma^-}$$
$$\{\sigma_{ij} v_{lj}\}|_{\Gamma^+} = -\{\sigma_{ij} v_{lj}\}|_{\Gamma^-}.$$

Then, applying spatial operator of (43) to the trace $\{\tilde{u}_i\}$ of a solution of (1) on grid (35) in each block, we have the following relationships for points at $\Gamma$ (after multiplication by $J^{-1}|T_{(l)}^{-1}|$):

$$\{-J^{-1}|T_{(l)}^{-1}|E_{ii'}\tilde{u}_{i'}\}|_{\Gamma^+} =$$
$$= \{-HJ^{-1}|T_{(l)}^{-1}|D_k^- JT_{kj}C_{ijj'i'}T_{k'j'}D_k^+\tilde{u}_{i'} + Q_k|T_{(l)}^{-1}|T_{kj}C_{ijj'i'}T_{k'j'}D_k^+\tilde{u}_{i'}\}|_{\Gamma^+} \approx$$
$$\approx \left(-H|T_{(l)}^{-1}|P\frac{d^2\{\tilde{u}_i\}}{dt^2} + HJ^{-1}|T_{(l)}^{-1}|\{F_i\}g(t) + \{\delta_{lk}Q_k|T_{(l)}^{-1}|T_{kj}C_{ijj'i'}T_{k'j'}D_k^+\tilde{u}_{i'}\}\right)\bigg|_{\Gamma^+}$$

and

$$\{-J^{-1} | T_{(l)}^{-1} | E_{ii'} \tilde{u}_{i'} \} |_{\Gamma^-} \approx$$
$$\approx \left( -H | T_{(l)}^{-1} | P \frac{d^2 \{\tilde{u}_i\}}{dt^2} + HJ^{-1} | T_{(l)}^{-1} | \{F_i\} g(t) + \{\delta_{lk} Q_k | T_{(l)}^{-1} | T_{kj} C_{ijj'i'} T_{k'j'} D_k^+ \tilde{u}_{i'} \} \right) \Big|_{\Gamma^-}.$$

The Kronecker symbol $\delta_{lk}$ is used here to emphasize that only terms with index $l$ remain at $\Gamma$.

Since $|T_{(l)}^{-1}| \{T_{lj} \sigma_{ij}\}|_{\Gamma^+} = |T_{(l)}^{-1}| \{T_{lj} \sigma_{ij}\}|_{\Gamma^-}$ and $Q_l|_{\Gamma^+} = -Q_l|_{\Gamma^-}$ then after summation of these expressions, taking into account annihilation of terms $\{\delta_{lk} Q_k | T_{(l)}^{-1} | T_{kj} C_{ijj'i'} T_{k'j'} D_k^+ \tilde{u}_{i'} \}$ at $\Gamma$, we have

$$\{-J^{-1} | T_{(l)}^{-1} | E_{ii'} \tilde{u}_{i'} \} |_{\Gamma^+} + \{-J^{-1} | T_{(l)}^{-1} | E_{ii'} \tilde{u}_{i'} \} |_{\Gamma^-} \approx$$
$$\approx \left( -H | T_{(l)}^{-1} | P \frac{d^2 \{\tilde{u}_i\}}{dt^2} + HJ^{-1} | T_{(l)}^{-1} | \{F_i\} g(t) \right) \Big|_{\Gamma^+}$$
$$+ \left( -H | T_{(l)}^{-1} | P \frac{d^2 \{\tilde{u}_i\}}{dt^2} + HJ^{-1} | T_{(l)}^{-1} | \{F_i\} g(t) \right) \Big|_{\Gamma^-}.$$

Therefore, we formulate the following equation in the grid points belonging $\Gamma$:

$$H | T_{(l)}^{-1} | P \frac{d^2 \{u_i\}}{dt^2} \Big|_{\Gamma^+} + H | T_{(l)}^{-1} | P \frac{d^2 \{u_i\}}{dt^2} \Big|_{\Gamma^-}$$
$$- \{J^{-1} | T_{(l)}^{-1} | E_{ii'} u_{i'} \} |_{\Gamma^+} - \{J^{-1} | T_{(l)}^{-1} | E_{ii'} u_{i'} \} |_{\Gamma^-} = \qquad (52)$$
$$= HJ^{-1} | T_{(l)}^{-1} | \{F_i\} g(t) |_{\Gamma^+} + HJ^{-1} | T_{(l)}^{-1} | \{F_i\} g(t) |_{\Gamma^-}.$$

This equation, in general, does not have a symmetric matrix of spatial operator because the multiplier $J^{-1} | T_{(l)}^{-1} |$ may be discontinuous at $\Gamma$. However, there is an important case where the symmetry is preserved. Consider the chain of equalities for points on $\Gamma$:

$$J^{-1} = \left( \det \frac{\partial(x_1, x_2, x_3)}{\partial(\xi_1, \xi_2, \xi_3)} \right)^{-1} = \det \frac{\partial(\xi_1, \xi_2, \xi_3)}{\partial(x_1, x_2, x_3)} = T_{lj} a_{lj} = T_{(l)} v_{lj} a_{lj},$$

where $a_{lj}$ denote the algebraic complements of matrix $J^{-1}$. We obtain that $J^{-1} T_{(l)}^{-1} = v_{lj} a_{lj}$. Suppose that the coordinate transformation in adjacent blocks is the same on tangent directions to $\Gamma$. Then $a_{lj}$ are continuous across the boundary, so it follows from $v_{lj} a_{lj}|_{\Gamma^+} = -v_{lj} a_{lj}|_{\Gamma^-}$ that $J^{-1} | T_{(l)}^{-1} |\big|_{\Gamma^+} = J^{-1} | T_{(l)}^{-1} |\big|_{\Gamma^-}$. Therefore, the factor $J^{-1} | T_{(l)}^{-1} |$ in (52) can be reduced, and we obtain at $\Gamma$:

$$(HJP|_{\Gamma^+} + HJP|_{\Gamma^-})\frac{d^2\{u_i\}}{dt^2} - \{E_{ii'}u_{i'}\}|_{\Gamma^+} - \{E_{ii'}u_{i'}\}|_{\Gamma^-} = \\ = H\{F_i\}g(t)|_{\Gamma^+} + H\{F_i\}g(t)|_{\Gamma^-} \quad . \tag{53}$$

The proposed equation (53) for handling points on the boundaries together with equation (42) for points inside blocks form a symmetric matrix of the spatial operator.

**Remark**. Although the notation for the operators in formulas (53), (42) is the same for adjacent blocks, the operators $E_{ii'}$ can be different, of course. The requirements are only 1) the common boundary grid points, 2) the common coordinate transformation in adjacent blocks for the tangent directions, and 3) the same triplets of $D_j^+$, $D_j^-$, and $H^j$ operators in adjacent blocks for the tangent directions ($j \neq l$).

## 4.5  *On the stability and accuracy of proposed difference schemes*

Evidently, it is possible to formulate elastic IBVP with different combinations of boundary conditions, similarly to proposed problems (46), (48), (50). For example, a typical formulation for seismic applications includes the free surface condition at the top boundary and nonreflecting conditions on all other boundaries. Also, one can use several blocks and impose the transmission conditions (53) for common boundaries. In all cases, we will obtain symmetric matrices of the spatial operator according to the proposed approach. Therefore, we can formulate the following theorem due to above theory and the lemma.

**Theorem 4**. *The explicit scheme (9) of the time integration of the problems (46), (48), (50), as well as of similarly formulated problems on the basis of (42) and combinations of boundary conditions (45), (4), (49) (with(51)), and (53) at different faces of a block computational domain is stable under the condition (10).*

The accuracy order of solutions will remain the same as the order of approximation due to a priori estimate (11). However, as noted in Section 3.4, it is much more important to investigate the accuracy properties numerically. Some examples of such test tasks are given in Section 5.

## *4.6 SEM scheme as a combination of SBP and multiblock approaches*

As is known, the computational SEM domain consists of hexahedral cells equipped by internal subgrids. Each such subgrid is generated by a nondegenerate map of the canonical cube with the rectangular GLL grid onto the cell. Therefore, the common face $\Gamma$ of two adjacent cells has the same transformation in the tangential directions. The transformations can be different in the normal direction only. Therefore, the ratio of Jacobians in adjacent cells on the common boundary $\Gamma$ is equal to the ratio of linear size of infinitesimal volume in the normal direction, i.e.

$$\frac{J\big|_{\Gamma^+}}{J\big|_{\Gamma^-}} = \frac{|T_{(l)}^{-1}|\big|_{\Gamma^+}}{|T_{(l)}^{-1}|\big|_{\Gamma^-}},$$

see Fig. **1**. This means that the factor $J^{-1}|T_{(l)}^{-1}|$ is continuous on boundaries of adjacent cells. Using (42) inside SEM cells for SBP operators $D^+ = D$, $D^- = D^\mathrm{T}$ described in Section 3.1.3 and (53) for grid points on the common boundaries, we obtain approximation of (34) with symmetric matrix of the spatial operator (similar averaging formulas to (53) for common grid points at edges and corners are easily derived). Thus one can treat considered SEM implementation as a finite-difference multiblock algorithm whose blocks are cells with specific SBP approximation on GLL grid.

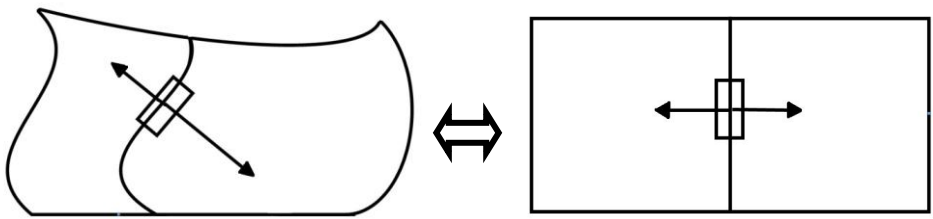

Fig. 1. One-to-one map of SEM adjacent cells to the canonical cells. Deformation of the boundary infinitesimal volume in the normal direction

# 5  Results of numerical experiments

The proposed difference schemes are implemented as a parallel MPI code on GPU cards for two and three dimensions. Below, we report three test results from numerous computations on studying stability and accuracy of the method.

## *5.1  A one-dimensional problem for the wave equation with a point source*

In this example, we consider the problem of spurious solutions and simultaneously demonstrate the advantages of the proposed nonsymmetric operators $D^+$, $D^-$.

We define the discrete delta function $\bar{\delta}_h(x_n)$ on the grid (18) taking its equal to the matrix element $H^{-1}$ at $x_n$ and to zero at the rest points.

Consider the wave problem on the interval $-1 \leq x \leq 1$ for equation (26) with a point source $\bar{\delta}_h(x_n)g(t)$, $x_n = 0$, in the RHS; the time impulse is generated by a Ricker wavelet $g(t)$. Initial conditions and Dirichlet boundary conditions are uniform. According to Section 3, we obtain the system

$$\mathrm{P}\bar{u}_{tt} - H^{-1}\tilde{L}\bar{u} = \bar{\delta}_h(x_n)g(t). \qquad (54)$$

Take $\rho = c = 1$ and consider two schemes. In the first scheme, we use the central differences $D$ of fourth-order approximation, i.e. $D^+ = D$, $D^- = D^\mathrm{T}$ here. In the second scheme, we use the operators (19) on the nonsymmetric stencils. Fig. **2** shows the solutions obtained by these two schemes. As we can see, the smooth part of the signal is the same, whereas the parasitic "saw" is visible for the central difference case only; it runs much faster than the main wave.

Obviously, it is the point source which is the reason for the observed large-amplitude parasitic waves in the first scheme as for a sufficiently smooth spatial source, all nonphysical solutions are within the approximation error. On the other hand, the proposed scheme with nonsymmetric operators $D^+$, $D^-$ does not generate such waves and has the necessary approximation properties.

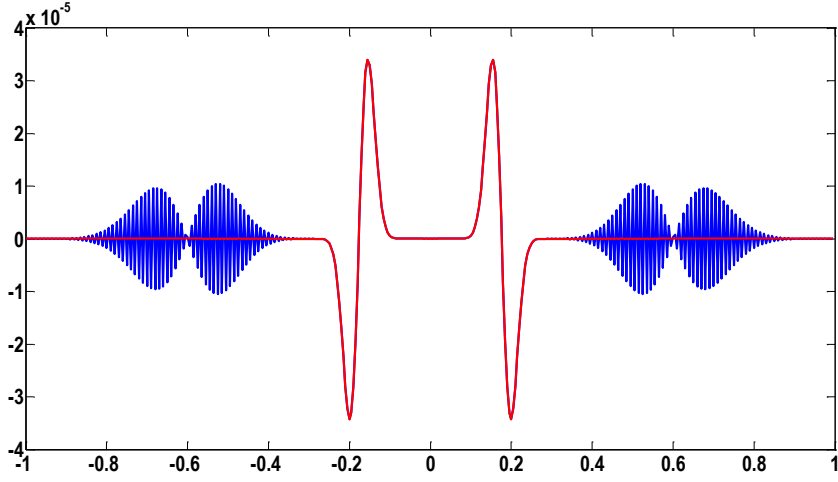

Fig. 2. Solution of the one-dimensional wave equation with a point source in space at $x=0$. Blue graph: central-difference operators ($D^+ = D$, $D^- = D^T$); red graph: nonsymmetric operators $D^+$, $D^-$

## *5.2 Lamb problem: example of two-dimensional elasticity*

Here we report the results of tests for the two-dimensional Lamb's problem [13]. The solution satisfies the Navier wave equation in the bottom half-plane with a homogeneous isotropic medium. A force point source is located at the surface with the non-zero component in the normal direction. The task is a good accuracy test because of the known analytical solution. We analyze the behavior of calculated solutions far and near to the point source.

Consider (1) in 2D formulation for $\boldsymbol{u} = (u_1, u_2)$ where the coordinate $x_2$ corresponds to the vertical axis. An elastic medium with constant density and velocities in the half-plane $x_2 \leq 0$ is described by Hooke's law with the stiffness tensor $C_{ijkl} = \lambda \delta_{ij}\delta_{kl} + \mu(\delta_{ik}\delta_{jl} + \delta_{il}\delta_{jk})$ where $\lambda = \rho(V_P^2 - 2V_S^2)$, $\mu = \rho V_S^2$. Free surface boundary conditions $\sigma_{ij}\nu_j = 0$ are imposed at the half-plane boundary $x_2 = 0$, where $(\nu_1, \nu_2)^T$ is the normal vector. Elastic waves are excited by the point source in the form of the normally directed force vector $\boldsymbol{F} = (\nu_2, -\nu_1)^T \delta(\boldsymbol{x} - \boldsymbol{x}_p) g(t)$ at $\boldsymbol{x}_p$. The analytic solution is known and described, e.g., in [7]. For calculating reference solutions, we use a free Fortran code available in [29] and based on approaches [2], [9].

Parameters of the isotropic medium are $V_P = 3200$ m/s, $V_S = 1847.5$ m/s, $\rho = 2200$ kg/m$^3$; the estimated Rayleigh wave velocity is $1698.6$ m/s. The time impulse is a Ricker wavelet $g(t) = (1 - 2\pi^2 \nu_0^2 t^2) \exp(-\pi^2 \nu_0^2 t^2)$ with the central frequency $\nu_0 = 10 Hz$. The computational domain has dimensions (–1500 to 1500)×(–1500 to 0) in meters; point source has coordinates (0,0). Boundary conditions (49) are imposed for open boundaries; actually they do not influence the solution at receivers for the considered simulation time. The basic grid has 251×181 nodes with 12 m and 8.3 m spacing in horizontal and vertical directions, respectively. The time integration step $\tau = 0.4$ ms. *C*-norm for the vector amplitude $\sqrt{u_1^2 + u_2^2}$ is used for accuracy measurements.

*5.2.1 Solution in the far field*

We compare the numerical solution with the analytical solution on the surface at the point with coordinates (600,0). For the fourth-order scheme on the boundary (the eighth-order in interior points) we observe the 1.2% relative accuracy. For the second-order scheme (fourth inside), the accuracy is significantly worse, approximately 20%. The vertical component $u_{h,2}(t,600,0)$ of the solution for these calculations is shown in Fig. **3**. Note that to achieve 1.2% accuracy by the second-order scheme, we need a 751 × 421 grid with 4 m and 3.6 m spatial spacing in the horizontal and vertical directions, respectively. Simulation on this grid requires 7 times more memory and 3.7 times more wall clock time.

We consider also the case of an inclined surface; that this is a serious problem for most finite difference algorithms using rectangular grids. An analytic solution for the inclined surface is obtained by rotating the coordinate system to the appropriate angle. To adapt our grid to this boundary, we use parametric coordinates. The oblique 251 × 181 grid is shown in Fig. **4**, each 10th grid line. Two snapshots of the solution amplitude are also presented here (fourth-order scheme). Relative accuracy of the solution at the point with coordinates (591, 104), i.e., on the surface at a distance 600 m from the source, is 1.3%. This slightly larger error is because of the small increase in grid spacing in the physical region with the oblique boundary.

We calculate also the solution for the oblique case on the finer grid with 501 × 361 points and with time step $\tau = 0.1$ ms to verify the convergence order. The resulting relative error is 0.06%, which is even better than the theoretical value estimated as 1.3 / 16 = 0.08%.

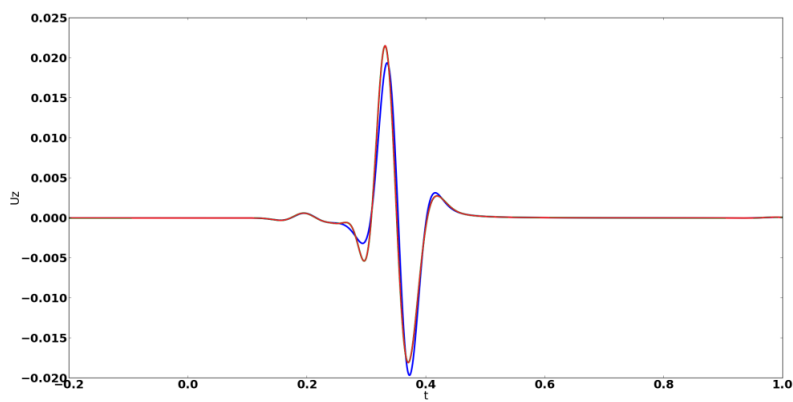

Fig. 3. Solution component $u_{h,2}(t,600,0)$. The grid has 251×181 points; the time step is 0.4 ms. The blue curve is for the fourth-order scheme on the boundary. The red curve is for the second-order scheme.

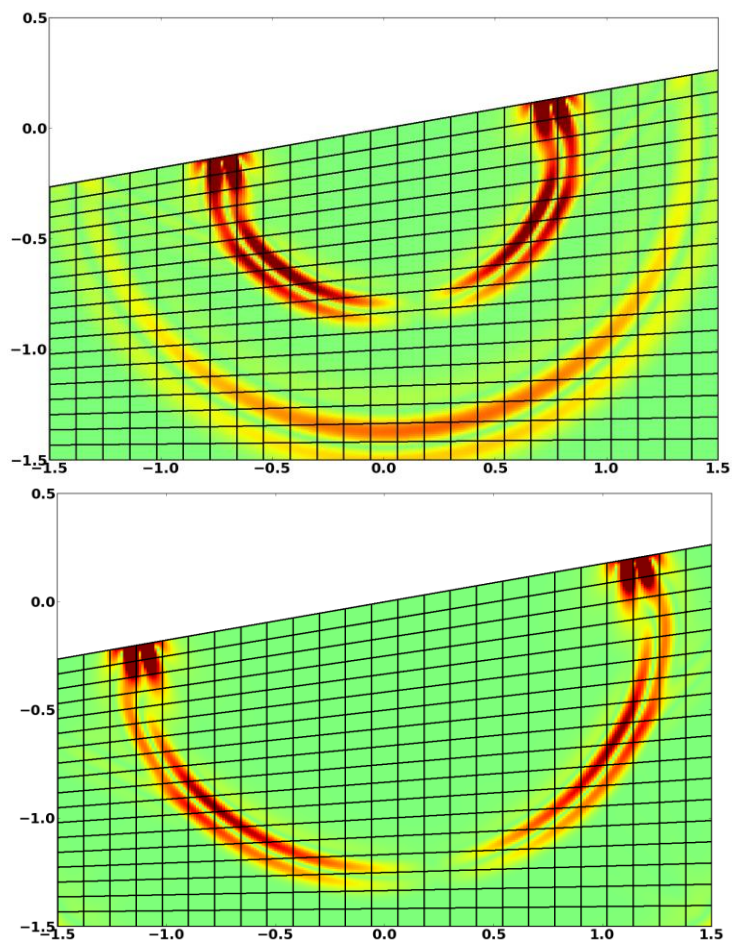

Fig. 4. Solution amplitude for the inclined surface case at the instants 0.44 s (left) and 0.68 s (right). Every 10[th] grid line is shown

*5.2.2 Accuracy of the calculation with the point source*

To analyze properties of the proposed difference scheme while simulating the point source, we examine the accuracy dependence versus the distance from the source. Fig. **5** shows the relative (percentage) error in C-norm for the case of a horizontal free surface depending on the distance from the source. The blue curve corresponds to the grid $251 \times 181$ with 12m horizontal spacing. We see that the error is less than 1% starting from 96 meters from the source. This distance is equal to $p = 8$ grid spacings, i.e., to the width of a standard central-difference stencil approximating the second derivative with $p^{th}$ order of accuracy in the horizontal direction. Recall that our scheme has the $p^{th}$ and $(p/2)^{th}$ orders of accuracy at the free surface in the horizontal and vertical directions, respectively. Calculation with the twice finer grid (red line) shows that the error drops $2^{p/2} = 16$ times at the same distance (96m). At the distance of $p$ grid spacings, i.e., 48 m, the error is 0.5%.

This and further analysis of the curves in Fig. **5** leads to the following observations:

- The scheme has $O(h)$ accuracy on the distance of $h$ to $ph$ from the point source, $p = 8$.
- The scheme comes to $O(h^{p/2})$ accuracy at the distance farther than $2ph$ from the point source.

In our opinion, these accuracy properties of the scheme are good enough for the considered point source, taking into account that this case is beyond the conventional theory of high-order accurate difference schemes for smooth input data.

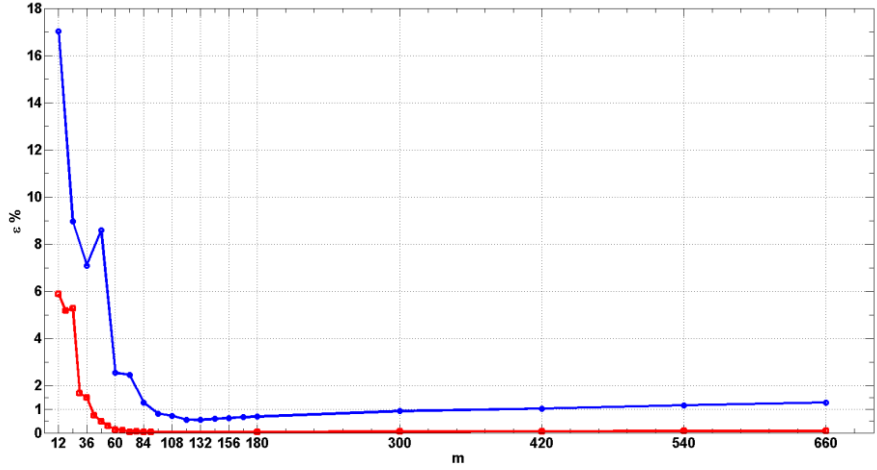

Fig. 5. Relative error of the solution amplitude on the horizontal surface versus the distance from a point source; calculations with grid spacing $h = 12$ m (blue line) and $h = 6$ m (red line)

## *5.3  3D modeling with anisotropy*

We demonstrate the main features of the proposed algorithm on an example of the seismic wave modeling in a two-layered anisotropic elastic medium with topography. Both horizontal surfaces are generated using Gaussian valleys and hills with random parameters, see Fig. **6** (land surface elevation amplitude is 1450 m); the bottom surface shape repeats the interface shape. A two-dimensional vertical cross-section of the computational grid is also shown in Fig. **6**. The horizontal dimension of the computational domain is 6 km × 6 km. We apply a smooth coordinate transformation to adapt the computational grid for three sub-horizontal boundaries and to reduce vertical spacing nearby them to compensate for the half-order accuracy of the scheme in this region. The grid has also a larger spacing near the open boundaries to save memory. Free surface boundary conditions are imposed on the top surface and non-reflecting conditions on the other boundary surfaces. Thomsen's parameters are $V_P = 2\,\text{km/s}$, $V_S = 1.2\,\text{km/s}$, $\rho = 2\,\text{g/cm}^3$, $\varepsilon = 0.334$, $\gamma = 0.575$, $\delta = 0.818$, $\alpha = \beta = 45°$ for the top TTI medium, and $V_P = 3\,\text{km/s}$, $V_S = 1.6\,\text{km/s}$, $\rho = 2\,\text{g/cm}^3$, $\varepsilon = 0.022$, $\gamma = 0.087$, $\delta = -0.072$, $\alpha = 90°$, $\beta = 15°$ for the bottom TTI medium. An explosion source with a 10 Hz central frequency Ricker wavelet is located in the center of the land surface. The grid volume is $501\times501\times(141+111)$ nodes; the scheme has the eighth-order in internal points and the fourth-order in the normal direction near boundaries. Snapshots of the solution amplitude are shown in Fig. 7. We clearly see the complex structure of different surface and body waves.

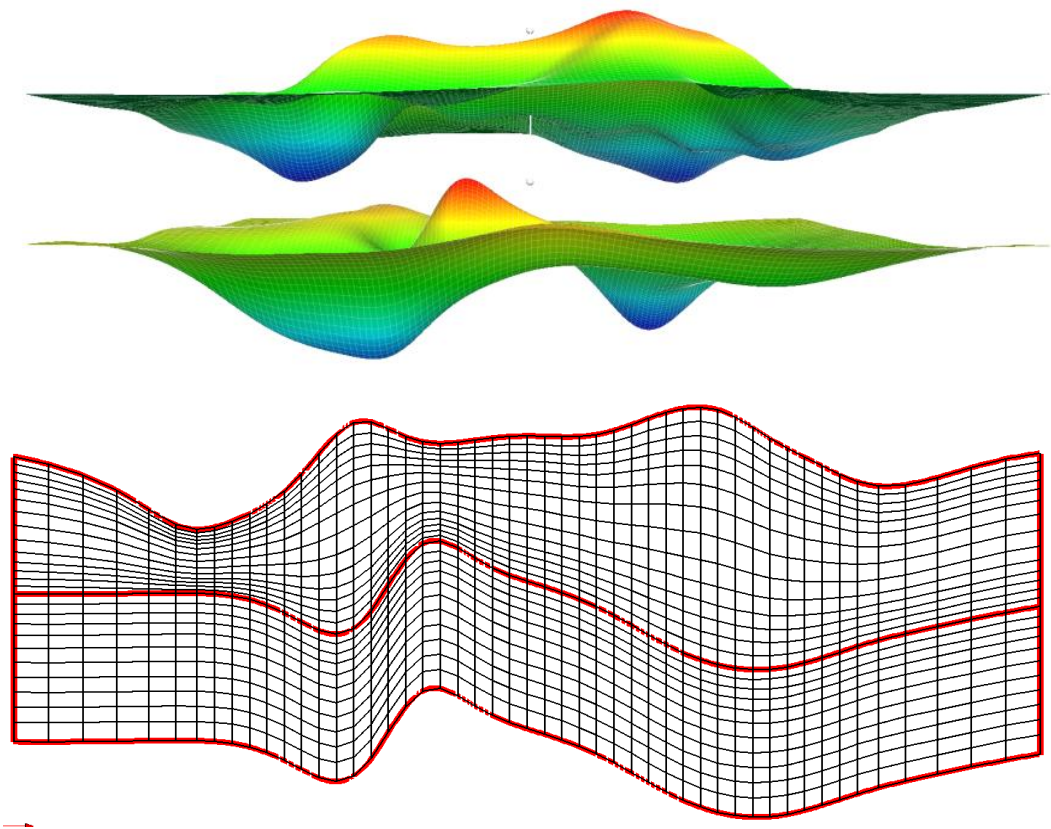

Fig. 6. Geometry of land surface and sub-horizontal interface (top). Vertical cross-section of the computational grid, each 10$^{th}$ line (bottom)

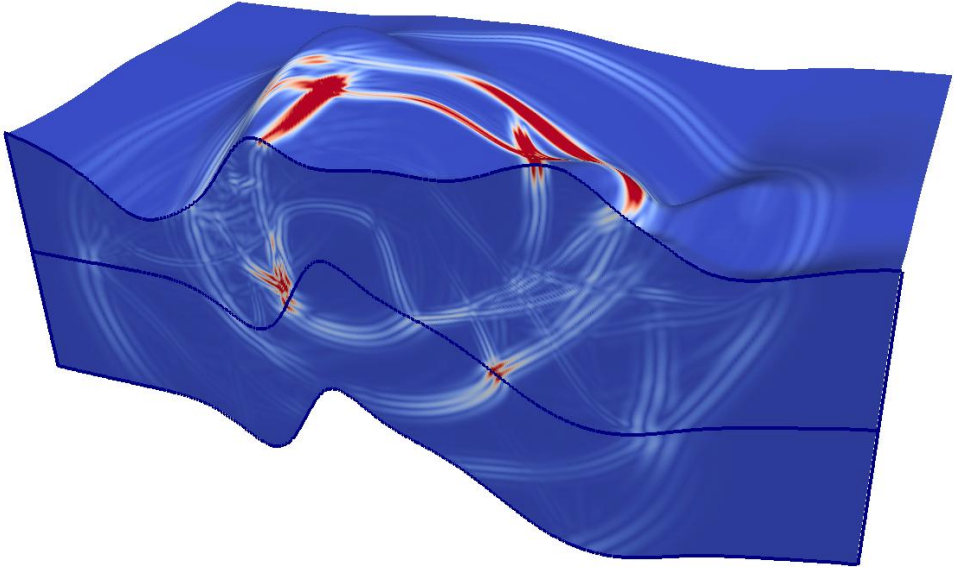

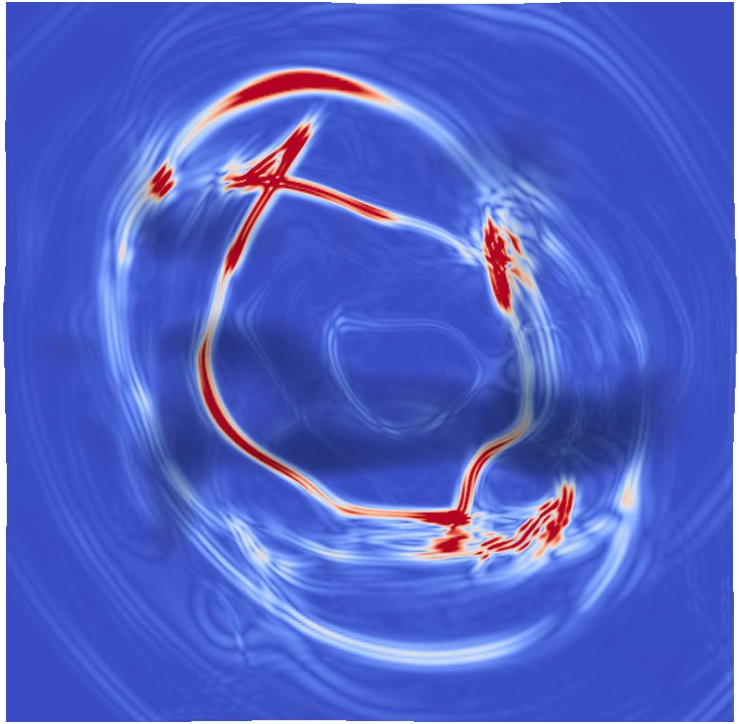

Fig. 7. Amplitude of displacement vector in the computational domain at $t=1.06$ s (top) and on the land surface at $t=1.56$ s, view from the top (bottom)

## 6   Conclusions

The paper presents the high-order accurate finite-difference scheme for solving three-dimensional problems of elastodynamics in anisotropic heterogeneous media. Approximation of the Navier wave equation is made in parametric coordinates, which allows obtaining curvilinear grids in physical coordinates adapted to geometry and velocity heterogeneities. The scheme also contains the approximation formulas of transmission conditions on the boundaries of the multiblock approach. We use a conventional explicit second-order central-difference time integration operator and obtain stable runs under the CFL condition.

The scheme has up to $p^{th}$ spatial accuracy order in the interior and $(p/2)^{th}$ order in normal direction to the boundaries, $p=4,6,8$, for uniform grids in parametric coordinates. We propose forward and backward finite differences derived in frames of the modified SBP rule to attenuate spurious saw-tooth solutions when solving governing equations with nonsmooth coefficients and/or point source RHS. Because of the use of sequential approximation of higher derivatives by the first-order difference operators, we optimize the theoretical number of operations per grid point, reduce the amount of memory in the implementation of the Hooke's law, and build an efficient parallel algorithm for multicore computing systems.

We also show that the conventional spectral element method belongs to class of derived schemes in the multiblock framework whose blocks are the SEM cells with SBP operators on GLL grid.

Numerous tests show the expected convergence rate to solutions on finer and finer grids as well as the long computation time stability, including multiblock cases. In particular, we describe the 2D Lamb problem modeling results to demonstrate good approximation properties of the method in the far and near field of a point source in space. Also, we give an illustration of the 3D solution in a heterogeneous anisotropic medium with a curvilinear free surface and multiblock grid.

# 7 Acknowledgments

The authors are grateful to Schlumberger for permission to publish the work. The first author is also grateful for support of the joint scientific-educational program between Schlumberger and the Moscow Institute of Physics and Technology.

# References

[1] D. Appelö, and N.A. Petersson, A stable finite difference method for the elastic wave equation on complex geometries with free surfaces. Commun. Comput. Phys., 5 (2009), pp. 84–107.

[2] L. Cagniard, Reflexion et refraction des ondes s´eismiques progressives, Paris, 1939

[3] M.H. Carpenter, D. Gottlieb, and S. Abarbanel, The stability of numerical boundary treatments for compact high-order finite-difference schemes. J. Comput. Phys., 108 (2) (1994), pp. 272-295

[4] L.E. Dovgilovich, and I.L. Sofronov, High-order finite-difference method of calculating wavefields in anisotropic media, Seismic Technologies, No.2 (2013), pp. 24-30 (in Russian).

[5] L. Dovgilovich, High-order FD method on curvilinear grids for anisotropic elastodynamic simulations, 75th EAGE Conference & Exhibition, Extended Abstracts (2013).

[6] G.J. Gassner, A skew-symmetric discontinuous Galerkin spectral element discretization and its relation to SBP-SAT finite difference methods. SIAM J. Sci. Comput. 35 (3), 2013, pp. A1233–A1253.

[7] A.G. Gorshkov, A.L. Medvedskii, L.N. Rabinskii, D.V. Tarlakovskii , Waves in solids, Moscow, "FizMatLit" (2004), pp. 1-472 (in Russian).

[8] O. Holberg, Computational aspects of the choice of operator and sampling interval for numerical differentiation in large-scale simulation of wave phenomena, Geophys. Prospect., 35 (1987), pp. 629–655.

[9] A.T. de Hoop, A modification of Cagniard's method for solving seismic pulse problems, Appl. Sci. Res. Sect. B 8 (1960), V. 4, pp. 349–356.

[10] D. Komatitsch, J.P. Vilotte, The spectral element method: an efficient tool to simulate the seismic response of 2D and 3D geological structures, Bull. Seis. Soc. Am., V. 88 (1998), pp. 368–392.

[11] J.E. Kozdon, E.M. Dunham, and J. Nordström, Simulation of dynamic earthquake ruptures in complex geometries using high-order finite difference methods, J. Scien. Comp. 55 (1), (2013), pp. 92-124.

[12] H.-O. Kreiss, and G. Scherer, Finite element and finite difference methods for hyperbolic partial differential equations, Mathematical Aspects of Finite

Elements in Partial Differential Equations. – Academic Press, (1974) pp. 195-212.

[13] H. Lamb, On the propagation of tremors over the surface of an elastic solid, Phi. Trans. Roy. Soc. (London), A 203 (1904), pp. l-42.

[14] S.G. Lekhnitskii, Theory of Elasticity of an anisotropic body, "Mir" (1981), pp. 1-431.

[15] V. Lisitsa, and D. Vishnevskiy, Lebedev scheme for the numerical simulation of wave propagation in 3D anisotropic elasticity. Geophysical Prospecting, 58(4) (2010), pp. 619–635.

[16] K. Mattsson, F. Ham, and G. Iaccarino, Stable boundary treatment for the wave equation on second-order form, J. Scien. Comp., 41(3), (2009), pp. 366–383.

[17] K. Mattsson, and J. Nordström, Summation by parts operators for finite difference approximations of second derivatives, J. Comp. Phys., 199 (2004), pp. 503-540.

[18] K. Mattsson, Summation by parts operators for finite difference approximations of second-derivatives with variable coefficients, J. Scien. Comp., 51(3), (2012), pp. 650–682.

[19] J. Nordström, K. Mattsson, and M. Svärd, Stable and accurate artificial dissipation, J. Scien. Comp., 21(1) (2004), pp. 57–79.

[20] N.A. Petersson, and B. Sjögreen, An energy absorbing far-field boundary condition for the elastic wave equation, Commun. Comput. Phys., 6 (2009), pp. 483–508.

[21] V.S. Ryaben'kii, The method of intrinsic boundary conditions in the theory of difference boundary value problems, Uspekhi Mat. Nauk, 26:3(159) (1971), pp. 105–160.

[22] V.S. Ryaben'kii, Method of difference potentials and its applications, Springer (2001) pp. 1- 556.

[23] E.H. Saenger, N. Gold, and S.A. Shapiro, Modeling the propagation of elastic waves using a modified finite-difference grid, Wave Motion, 31 (2000), pp. 77-92.

[24] A.A. Samarski, The theory of difference schemes. CRC Press, (2001), pp. 1-786.

[25] I.L. Sofronov, Differential part of transparent boundary conditions for certain hyperbolic systems of second-order equations, Doklady Mathematics, 79(3) (2009), pp. 412–414.

[26] I.L. Sofronov, N.A. Zaitsev, A. Daryin, O. Voskoboinikova, and L. Dovgilovich, Multi-block FD method for 3D geophysical simulation with explicit representation of sub-horizontal interfaces. 74th EAGE Conference & Exhibition, Extended Abstracts (2012).

[27] B. Strand, Summation by parts for finite difference approximations for d/dx. J. Comput. Physics, 110 (1994), pp. 47–67.

[28] J. Virieux, P-SV wave propagation in heterogeneous media: velocity-stress finite-difference method, Geophysics, 51 (1986), pp. 889-901.

[29] http://www.spice-rtn.org/library/software/softwarefolder_view.html  (2013)